%&amstex
\input amstex
\input amsppt.sty
\magnification=\magstep1
\hsize=30truecc
\vsize=22.2truecm
\baselineskip=16truept
\NoBlackBoxes
\TagsOnRight \pageno=1 \nologo
\def\Z{\Bbb Z}
\def\N{\Bbb N}

\def\l{\left}
\def\r{\right}
\def\bg{\bigg}
\def\({\bg(}
\def\[{\bg\lfloor}
\def\){\bg)}
\def\]{\bg\rfloor}
\def\t{\text}
\def\f{\frac}

\def\bi{\binom}
\def\eq{\equiv}

\def\ls{\leqslant}
\def\gs{\geqslant}
\def\mo{\roman{mod}}

\def\al{\alpha}
\def\da{\delta}

\def\Proof{\noindent{\it Proof}}

\def\Remark{\medskip\noindent{\it  Remark}}

\hbox {Preprint, {\tt arXiv:0910.5667}}
\bigskip
\topmatter
\title On sums of binomial coefficients modulo $p^2$\endtitle
\author Zhi-Wei Sun\endauthor
\leftheadtext{Zhi-Wei Sun} \rightheadtext{On sums of binomial coefficients modulo $p^2$}
\affil Department of Mathematics, Nanjing University\\
 Nanjing 210093, People's Republic of China
  \\  zwsun\@nju.edu.cn
  \\ {\tt http://math.nju.edu.cn/$\sim$zwsun}
\endaffil
\abstract Let $p$ be an odd prime and $a$ be a positive integer. In this paper
we investigate the sum $\sum_{k=0}^{p^a-1}\bi{hp^a-1}k\bi{2k}k/m^k$ mod $p^2$,
where $h,m$ are $p$-adic integers with $m\not\eq0\ (\mo\ p)$. For example, we show that if
$h\not\eq0\ (\mo\ p)$ and $p^a>3$, then
 $$\align&\sum_{k=0}^{p^a-1}\bi{hp^a-1}k\bi{2k}k\l(-\f h2\r)^k
\\\eq&\l(\f{1-2h}{p^a}\r)\(1+h\(\l(4-\f 2h\r)^{p-1}-1\)\)\ (\mo\ p^2),
\endalign$$
where $(-)$ denotes the Jacobi symbol.
 Here is another remarkable congruence:
If $p^a>3$ then
$$\sum_{k=0}^{p^a-1}\bi{p^a-1}k\bi{2k}k(-1)^k\eq 3^{p-1}\l(\f{p^a}3\r)\ (\mo\ p^2).$$
\endabstract
\thanks 2010 {\it Mathematics Subject Classification}.\,Primary 11B65;
Secondary 05A10,\,11A07, 11B39, 11S99.
\newline\indent Supported by the National Natural Science
Foundation (grant 10871087) and the Overseas Cooperation Fund (grant 10928101) of China.
\endthanks
\endtopmatter
\document

\heading{1. Introduction}\endheading

Let $p$ be a prime.
In 2006  H. Pan and Z. W. Sun [PS]
proved that
$$\sum_{k=0}^{p-1}\bi{2k}{k+d}\eq\l(\f{p-d}3\r)\ (\mo\ p)\quad \t{for}\ d=0,\ldots,p,$$
where $(-)$ is the Jacobi symbol. Later Sun and R. Tauraso [ST1, ST2] determined
$\sum_{k=0}^{p^a-1}\bi{2k}k/m^k$ and $\sum_{k=1}^{p-1}\bi{2k}k/(km^{k-1})$ modulo $p$
via Lucas sequences, where $m$ is an integer not divisible by $p$ and $a$ is any positive integer.
Recently Sun [S10a] was able to determine
$\sum_{k=0}^{p^a-1}\bi{2k}k/m^k$ mod $p^2$. See also [SSZ], [GZ] and [S10b] for related results on $p$-adic valuations.

In this paper we study a new kind of sums
$$\sum_{k=0}^{p^a-1}\bi{hp^a-1}k\f{\bi{2k}k}{m^k}$$
modulo $p^2$, where $a\in\Z^+=\{1,2,3,\ldots\}$, and $h$ and $m$ are $p$-adic integers with $m\not\eq0\ (\mo\ p)$.

For a prime $p$ we use $\Z_p$ to denote the ring of $p$-adic integers;
if $h\in\Z_p$ and $h\not\eq0\ (\mo\ p)$ then we denote the quotient $(h^{p-1}-1)/p\in\Z_p$
by $q_p(h)$ and call it a {\it Fermat quotient}. For $m,n\in\N=\{0,1,2,\ldots\}$, the Kronecker symbol $\da_{m,n}$
takes $1$ or $0$ according as $m=n$ or not.

Now we state our main results and give some corollaries.

\proclaim{Theorem 1.1}
Let $p$ be an odd prime and let $a\in\Z^+$. Let $h$ be a $p$-adic integer with $h\not\eq0\ (\mo\ p)$,
and $(2h\not\eq1\ (\mo\ p)$ or $p^a>3)$. Then
$$\aligned&\sum_{k=0}^{p^a-1}\bi{hp^a-1}k\bi{2k}k\l(-\f h2\r)^k
\\\eq&\l(\f{1-2h}{p^a}\r)\(1+h\(\l(4-\f2h\r)^{p-1}-1\)\)\ (\mo\ p^2).
\endaligned\tag1.1$$
\endproclaim

\proclaim{Corollary 1.1} Let $p$ be an odd prime and let $a\in\Z^+$. Then
$$\sum_{k=0}^{p^a-1}\bi{p^a-1}k\f{\bi{2k}k}{(-2)^k}\eq(-1)^{(p^a-1)/2}2^{p-1}\ (\mo\ p^2).\tag1.2$$
\endproclaim
\Proof. Simply apply Theorem 1.1 with $h=1$. \qed

\Remark\ 1.1. Let $m\in\Z$ and $n\in\Z^+$. Later we will show that
$$\aligned&\sum_{k=0}^{n-1}\bi{n-1}k\bi{2k}k(-1)^km^{n-1-k}
\\=&\sum_{k=0}^{\lfloor(n-1)/2\rfloor}\bi{n-1}k\bi{n-1-k}k(m-2)^{n-1-2k}.
\endaligned\tag1.3$$
Thus, for any prime $p>3$, by applying Morley's congruence (cf. [M], [C] and [P])
$$\bi{p-1}{(p-1)/2}\eq (-1)^{(p-1)/2}4^{p-1}\ (\mo\ p^3)$$
we get
$$\sum_{k=0}^{p-1}\bi{p-1}k\f{\bi{2k}k}{(-2)^k}
\eq(-1)^{(p-1)/2}2^{p-1}\ (\mo\ p^3)$$
which is a refinement of (1.2) in the case $a=1$.

\proclaim{Corollary 1.2} Let $p>3$ be a prime and let $a\in\Z^+$. Then
$$\sum_{k=0}^{p^a-1}\bi{2p^a-1}k\bi{2k}k(-1)^k\eq\l(\f{p^a}3\r)(2\times 3^{p-1}-1)\ (\mo\ p^2)\tag1.4$$
and
$$\sum_{k=0}^{p^a-1}\bi{p^a+k}k\f{\bi{2k}k}{2^k}\eq\l(\f{3}{p^a}\r)(1-p(q_p(2)+q_p(3))\ (\mo\ p^2).\tag1.5$$
\endproclaim
\Proof. Just put $h=2$ and $h=-1$ in (1.1) and note that $\bi{-x}k=(-1)^k\bi{x+k-1}k$. \qed

\proclaim{Corollary 1.3} Let $p$ be an odd prime and let $a\in\Z^+$. Then
$$\sum_{k=0}^{p^a-1}\bi{2p^a+k}k\bi{2k}k(-1)^k\eq\l(\f{p^a}5\r)\l(3-2\times 5^{p-1}\r)\ (\mo\ p^2).\tag1.6$$
\endproclaim
\Proof. Simply apply (1.1) with $h=-2$. \qed

Now we need to introduce Lucas sequences.

Let $A,B\in\Z$. The Lucas sequences $u_n=u_n(A,B)\ (n\in\N)$ and  $v_n=v_n(A,B)\ (n\in\N)$
are defined by
$$u_0=0,\ u_1=1,\ \t{and}\ u_{n+1}=Au_n-Bu_{n-1}\ (n=1,2,3,\ldots)$$
and
$$v_0=2,\ v_1=A,\ \t{and}\ v_{n+1}=Av_n-Bv_{n-1}\ (n=1,2,3,\ldots).$$
The characteristic equation $x^2-Ax+B=0$ has two roots
$$\al=\f{A+\sqrt{\Delta}}2\quad\t{and}\quad\beta=\f{A+\sqrt{\Delta}}2,$$
where $\Delta=A^2-4B$. It is well known that for any $n\in\N$ we have
$$u_n=\sum_{0\ls k<n}\al^k\beta^{n-1-k}\quad\t{and }\quad v_n=\al^n+\beta^n.$$
If $p$ is a prime then
$$v_p=\al^p+\beta^p\eq(\al+\beta)^p=A^p\eq A\ (\mo\ p).$$
It is also known that
$$u_p\eq\l(\f{\Delta}p\r)\ (\mo\ p)\ \ \t{and}\ \ u_{p-(\f{\Delta}p)}\eq0\ (\mo\ p)$$
for any prime $p$ not dividing $2B$. (See, e.g., [S10a, Lemma 2.3].)
The reader may consult [S06] for connections between Lucas sequences and quadratic fields.

 Our following result is more general than Theorem 1.1.

\proclaim{Theorem 1.2} Let $p$ be an odd prime and $m\in\Z$ with $p\nmid m$. Set $\Delta=m(m-4)$
and let $h\in\Z_p$. Then we have
$$\aligned&\sum_{k=0}^{p^a-1}\bi{hp^a-1}k\f{\bi{2k}k}{(-m)^k}
\\\eq&\l(\f{\Delta}{p^{a-1}}\r)\l(1-\f{hm}2\r)u_{p-(\f{\Delta}p)}(m-2,1)
\\&+\l(\f{\Delta}{p^a}\r)(1+h((m-4)^{p-1}-1))
\\&-\cases h(m-4)\ (\mo\ p^2)&\t{if}\ p^a=3\ \t{and}\ 3\mid m-1,
\\0\ (\mo\ p^2)&\t{otherwise}.
\endcases
\endaligned\tag1.7$$
In particular, if $hm\eq2\ (\mo\ p)$ then
$$\aligned&\sum_{k=0}^{p^a-1}\bi{hp^a-1}k\f{\bi{2k}k}{(-m)^k}
\\\eq&\l(\f{\Delta}{p^a}\r)(1+h((m-4)^{p-1}-1))
\\&+\cases m-4\ \ (\mo\ p^2)&\t{if}\ p^a=3\ \t{and}\ 3\mid m-1,\\0\ (\mo\ p^2)&\t{otherwise}.\endcases.
\endaligned\tag1.8$$
\endproclaim

\proclaim{Corollary 1.4} Let $p$ be an odd prime and let $a\in\Z^+$ with $p^a>3$. Then
$$\sum_{k=0}^{p^a-1}\bi{p^a-1}k\bi{2k}k(-1)^k\eq 3^{p-1}\l(\f{p^a}3\r)\ (\mo\ p^2)\tag1.9$$
and
$$\sum_{k=0}^{p^a-1}\bi{p^a-1}k\f{\bi{2k}k}{(-3)^k}\eq\l(\f{p^a}3\r)\ (\mo\ p^2).\tag1.10$$
\endproclaim
\Proof. Just apply (1.7) with $h=1$ and $m\in\{1,3\}$
and note  that $(-1)^{n-1}u_n(1,1)=u_n(-1,1)=(\f n3)$ for $n\in\N$. \qed

\proclaim{Corollary 1.5} Let $p\not=2,5$ be a prime and let $a\in\Z^+$. Then
$$\sum_{k=0}^{p^a-1}\bi{p^a-1}k\bi{2k}k\eq \l(\f{p^a}5\r)\l(5^{p-1}-3F_{p-(\f p5)}\r)\ (\mo\ p^2)\tag1.11$$
and
$$\sum_{k=0}^{p^a-1}\bi{p^a-1}k\f{\bi{2k}k}{(-5)^k}\eq\l(\f{p^a}5\r)\l(1-3F_{p-(\f p5)}\r)\ (\mo\ p^2),\tag1.12$$
where $\{F_n\}_{n\gs0}$ is the well-known Fibonacci sequence defined by
$$F_0=0,\ F_1=1,\ \t{and}\ F_{n+1}=F_n+F_{n-1}\ (n=1,2,3,\ldots).$$
\endproclaim
\Proof. Observe that $$(-1)^{n-1}u_n(-3,1)=u_n(3,1)=F_{2n}=F_nL_n,$$
where $L_n=v_n(1,-1)$. By [SS, Corollary 1] (or the proof of Corollary 1.3 of [ST1]), if $p\not=2,5$ then
$L_{p-(\f p5)}\eq 2\l(\f p5\r)\ (\mo\ p^2).$  In view of this,
if we apply (1.7) with $h=1$ and $m\in\{-1,5\}$ the we obtain the desired result. \qed

 To conclude this section we raise four conjectures.

\proclaim{Conjecture 1.1} Let $p$ be an odd prime and let $h$ be an integer with $2h-1\eq0\ (\mo\ p)$.
If $a\in\Z^+$ and $p^a>3$, then
$$\sum_{k=0}^{p^a-1}\bi{hp^a-1}k\bi{2k}k\l(-\f h2\r)^k\eq0\ (\mo\ p^{a+1}).$$
Also, for any $n\in\Z^+$ we have
$$\f1n\sum_{k=0}^{n-1}\bi{hn-1}k\bi{2k}k\l(-\f h2\r)^k\in\Z_p.$$
\endproclaim

\proclaim{Conjecture 1.2} Let $p>3$ be a prime. Then
$$\sum_{k=0}^{p-1}\bi{p-1}k\bi{2k}k((-1)^k-(-3)^{-k})\eq\l(\f p3\r)(3^{p-1}-1)\ (\mo\ p^3).$$
\endproclaim

\proclaim{Conjecture 1.3} Let $p$ be a prime with $p\eq\pm1\ (\mo\ 12)$. Then
$$\sum_{k=0}^{p-1}\bi{p-1}k\bi{2k}k(-1)^ku_k(4,1)\eq(-1)^{(p-1)/2}u_{p-1}(4,1)\ (\mo\ p^3).$$
\endproclaim

\proclaim{Conjecture 1.4} Let $p$ be a prime with $p\eq\pm1\ (\mo\ 8)$. Then
$$\sum_{k=0}^{p-1}\bi{p-1}k\bi{2k}k\f{u_k(4,2)}{(-2)^k}\eq(-1)^{(p-1)/2}u_{p-1}(4,2)\ (\mo\ p^3).$$
\endproclaim

We remark that the author could prove the congruences in Conjectures 1.2--1.4 modulo $p^2$.

\medskip

In the next section we provide some lemmas. Section 3 is devoted to our proofs of Theorems 1.1--1.2 and (1.3).

\heading{2. Some Lemmas}\endheading

For $n\in\N$ we set $H_n=\sum_{0<k\ls n}1/k$.

\proclaim{Lemma 2.1} Let $p$ be an odd prime and let $m\in\Z$ with $p\nmid m$.
If $p\mid m-4$ then
$$\sum_{k=1}^{p^a-1}\f{p^{a-1}H_k}{m^k}\bi{2k}k\eq 2\da_{a,1}\ (\mo\ p).\tag2.1$$
When $m\not\eq4\ (\mo\ p)$, we have
$$\sum_{k=1}^{p^a-1}\f{p^{a-1}H_k}{m^k}\bi{2k}k\eq-\l(\f{m(m-4)}{p^a}\r)
\sum_{k=1}^{p-1}\f{\bi{2k}k}{k(4-m)^k}\ (\mo\ p).\tag2.2$$
\endproclaim
\Proof. For $k=1,\ldots,(p^a-1)/2$, we have
$$\align\f{\bi{(p^a-1)/2}k}{\bi{2k}k/(-4)^k}=&\f{\bi{(p^a-1)/2}k}{\bi{-1/2}k}=\prod_{j=1}^k\f{(p^a-1)/2-j+1}{-1/2-j+1}
\\=&\prod_{j=1}^k\l(1-\f{p^a}{2j-1}\r)\eq 1\ (\mo\ p).
\endalign$$
If $k\in\{(p^a+1)/2,\ldots,p^a-1\}$, then
$2k-p^a\in\{1,\ldots,k-1\}$ and hence
$$\bi{2k}k=\bi{p^a+(2k-p^a}{k}\eq\bi{p^a}0\bi{2k-p^a}k=0\ (\mo\ p)$$
with the help of Lucas' congruence. So, for any $k=0,\ldots,p^a-1$
we have
$$\bi{2k}k\eq(-4)^k\bi{(p^a-1)/2}k\ (\mo\ p).\tag2.3$$
Note also that
$$p^{a-1}H_k=\sum_{j=1}^k\f{p^{a-1}}j\in\Z_p\quad\t{for every}\
k=1,\ldots,p^a-1.$$
Therefore
$$\sum_{k=1}^{p^a-1}\f{p^{a-1}H_k}{m^k}\bi{2k}k\eq\sum_{k=1}^{(p^a-1)/2}\bi{(p^a-1)/2}k\l(-\f 4m\r)^k(p^{a-1}H_k)\ (\mo\ p).$$

For each $k\in\N$ clearly
$$\align H_k=&\sum_{0<j\ls k}\int_0^1x^{j-1}dx=\int_0^1\sum_{0<j\ls k}x^{j-1}dx
\\=&\int_0^1\f{1-x^k}{1-x}dx=\int_0^1\f{1-(1-t)^k}tdt.
\endalign$$
Thus
$$\sum_{k=1}^{p^a-1}\f{p^{a-1}H_k}{m^k}\bi{2k}k
\eq p^{a-1}\Sigma\ (\mo\ p),$$
where
$$\align \Sigma:=&\int_0^1\sum_{k=0}^{(p^a-1)/2}\bi{(p^a-1)/2}k\l(-\f 4m\r)^k\f{1-(1-t)^k}tdt
\\=&\int_0^1\f{(1-4/m)^{(p^a-1)/2}-(1-(1-t)4/m)^{(p^a-1)/2}}tdt
\\=&-\sum_{k=1}^{(p^a-1)/2}\bi{(p^a-1)/2}k\l(1-\f 4m\r)^{(p^a-1)/2-k}\int_0^1\l(\f{4t}m\r)^k\f{dt}t
\\=&-\f1{m^{(p^a-1)/2}}\sum_{k=1}^{(p^a-1)/2}\bi{(p^a-1)/2}k\f{4^k}{k}(m-4)^{(p^a-1)/2-k}.
\endalign$$

If $m\eq 4\ (\mo\ p)$, then
$$p^{a-1}\Sigma\eq -\f1{m^{(p^a-1)/2}}\cdot\f{p^{a-1}}{(p^a-1)/2}4^{(p^a-1)/2}\eq 2\da_{a,1}\ (\mo\ p)$$
and hence (2.1) holds.

Now assume that $m\not\eq 4\ (\mo\ p)$. In view of (2.3),
$$\align p^{a-1}\Sigma\eq&-\f{(m(m-4))^{(p^a-1)/2}}{m^{p^a-1}}\sum_{k=1}^{p^a-1}\bi{2k}k\f{(-1)^kp^{a-1}}{k(m-4)^k}
\\\eq&-\l(\f{m(m-4)}{p^a}\r)p^{a-1}\sum_{k=1}^{p^a-1}\f{\bi{2k}k}{k(4-m)^k}\ (\mo\ p).
\endalign$$
So it suffices to prove that
$$p^{a-1}\sum_{k=1}^{p^a-1}\f{\bi{2k}k}{kn^k}\eq \sum_{k=1}^{p-1}\f{\bi{2k}k}{kn^k}\ (\mo\ p)$$
for any $n\in\Z$ with $p\nmid n$. If $p^{a-1}\nmid k$ then $p^{a-1}/k\eq0\ (\mo\ p)$. Therefore
$$p^{a-1}\sum_{k=1}^{p^a-1}\f{\bi{2k}k}{kn^k}
\eq p^{a-1}\sum_{j=1}^{p-1}\f{\bi{2p^{a-1}j}{p^{a-1}j}}{p^{a-1}jn^{p^{a-1}j}}
\eq\sum_{j=1}^{p-1}\f{\bi{2j}j}{jn^j}\ (\mo\ p)$$
in view of the Lucas congruence (cf. [St, p.\,44]).

So far we have completed the proof of Lemma 2.1. \qed

\proclaim{Lemma 2.2 {\rm (Sun [S10a])}} Let $p$ be an odd prime and let $a\in\Z^+$. Let $m$ be any integer not divisible by $p$
and set $\Delta=m(m-4)$.
Then we have
$$\sum_{k=0}^{p^a-1}\f{\bi{2k}k}{m^k}\eq\l(\f{\Delta}{p^a}\r)
+\l(\f{\Delta}{p^{a-1}}\r)u_{p-(\f{\Delta}{p})}(m-2,1)\ (\mo\ p^2).$$
\endproclaim

\proclaim{Lemma 2.3 {\rm (Sun and Tauraso [ST1])}} Let $p$ be a prime and let $m$ be an
integer not divisible by $p$. Then we have
$$\f12\sum_{k=1}^{p-1}(-1)^k\f{\bi{2k}k}{km^{k-1}}\eq\f{m^p-v_p(m,-m)}p\
(\mo\ p).$$
\endproclaim

\proclaim{Lemma 2.4} Let $p$ be an odd prime and let $m\in\Z$ with $\Delta=m(m-4)\not\eq0\ (\mo\ p)$.
Then
$$\aligned&\f2{m-4}\cdot\f{v_p(m-4,4-m)-(m-4)^p}p
\\\eq&\f m2\l(\f{\Delta}p\r)\f{u_{p-(\f{\Delta}p)}(m-2,1)}p-q_p(m-4)\ (\mo\ p).
\endaligned\tag2.4$$
\endproclaim
\Proof. (i) Let us first show the equality
 $$\f{v_{2n+1}(m-4,4-m)}{(m-4)^{n+1}}=\f{u_{2n+1}(m,m)}{m^n}\tag2.5$$
 for $n=0,1,2,\ldots$. Clearly both sides of (2.3) coincide with $1$ when $n=0$.
 Note that
 $$\align &\f{v_3(m-4,4-m)}{(m-4)^2}=\f{v_2(m-4,4-m)+v_1(m-4,4-m)}{m-4}
 \\= &v_1(m-4),4-m)+v_0(m-4,4-m)+\f{v_1(m-4,4-m)}{m-4}
 \\=&m-4+2+1=m-1=u_2(m,m)-u_1(m,m)=\f{u_3(m,m)}{m}.
 \endalign$$
Also, for $n=2,3,\ldots$ we have
$$\align&\f{v_{2n+1}(m-4,4-m)}{(m-4)^{n+1}}
\\=&\f{v_{2n-1}(m-4,4-m)+v_{2n}(m-4,4-m)}{(m-4)^n}
\\=&\f{(1+(m-4))v_{2n-1}(m-4,4-m)+(m-4)v_{2n-2}(m-4,4-m)}{(m-4)^n}
\\=&\f{(m-2)v_{2n-1}(m-4,4-m)-(m-4)v_{2n-3}(m-4,4-m)}{(m-4)^n}
\\=&(m-2)\f{v_{2n-1}(m-4,4-m)}{(m-4)^n}-\f{v_{2n-3}(m-4,4-m)}{(m-4)^{n-1}}
\endalign$$
and
$$\align\f{u_{2n+1}(m,m)}{m^n}=&\f{u_{2n}(m,m)-u_{2n-1}}{m^{n-1}}
\\=&\f{(m-1)u_{2n-1}(m,m)-mu_{2n-2}(m,m)}{m^{n-1}}
\\=&\f{(m-1)u_{2n-1}(m,m)-(u_{2n-1}(m,m)+mu_{2n-3}(m,m))}{m^{n-1}}
\\=&(m-2)\f{u_{2n-1}(m,m)}{m^{n-1}}-\f{u_{2n-3}(m,m)}{m^{n-2}}.
\endalign$$
Thus, by induction (2.5) holds for all $n\in\N$.

(ii) By part (i),
$$u_p(m,m)=\f{m^{(p-1)/2}}{(m-4)^{(p+1)/2}}(v_p(m-4,4-m)-(m-4)^p)+(m(m-4))^{(p-1)/2}.
$$
Since $v_p(m-4,4-m)\eq (m-4)^p\ (\mo\ p)$ and
$$\align&\Delta^{(p-1)/2}-\l(\f{\Delta}p\r)
\\=&(m-4)^{(p-1)/2}\l(m^{(p-1)/2}-\l(\f mp\r)\r)+\l(\f mp\r)\((m-4)^{(p-1)/2}-\l(\f{m-4}p\r)\)
\\\eq&\l(\f{\Delta}p\r)\l(\f mp\r)\l(m^{(p-1)/2}-\l(\f mp\r)\r)
\\&+\l(\f{\Delta}p\r)\l(\f{m-4}p\r)\((m-4)^{(p-1)/2}-\l(\f{m-4}p\r)\)
\\\eq&\f12\l(\f{\Delta}p\r)\l(m^{p-1}-1+(m-4)^{p-1}-1\r)\ (\mo\ p^2),
\endalign$$
we have
$$\align u_p(m,m)-\l(\f{\Delta}p\r)\eq&\f{(\f mp)}{(m-4)(\f{m-4}p)}\l(v_p(m-4,4-m)-(m-4)^p\r)
\\&+\f12\l(\f{\Delta}p\r)\l(m^{p-1}-1+(m-4)^{p-1}-1\r)
\\\eq&\f1{m-4}\l(\f{\Delta}p\r)\l(v_p(m-4,4-m)-(m-4)^p\r)
\\&+\f p2\l(\f{\Delta}p\r)(q_p(m)+q_p(m-4))\ (\mo\ p^2).
\endalign$$
On the other hand, by [S10a, Lemma 2.4] we have
$$2u_p(m,m)-\l(\f{\Delta}p\r)m^{p-1}\eq u_p(m-2,1)+u_{p-(\f{\Delta}p)}(m-2,1)\ (\mo\ p^2).$$
Thus
$$\align &\f2{m-4}\l(\f{\Delta}p\r)(v_p(m-4,4-m)-(m-4)^p)
\\\eq&u_p(m-2,1)-\l(\f{\Delta}p\r)+u_{p-(\f{\Delta}p}(m-2,1)-\l(\f{\Delta}p\r)pq_p(m-4)\ (\mo\ p^2).
\endalign$$
Comparing this with (2.4) we have reduced (2.4) to the following congruence
$$u_p(m-2,1)-\l(\f{\Delta}p\r)\eq\l(\f m2-1\r)u_{p-(\f{\Delta}p)}(m-2,1)\ (\mo\ p^2).\tag2.6$$

Let $\al$ and $\beta$ be the two roots of the equation $x^2-(m-2)x+1=0$. Then
$$v_n(m-2,1)^2-\Delta u_n^2(m-2,1)=(\al^n+\beta^n)^2-(\al^n-\beta^n)^2=4(\al\beta)^n=4$$
for all $n\in\N$. As $u_{p-(\f{\Delta}p)}(m-2,1)\eq0\ (\mo\ p)$ we have
$$v_{p-(\f{\Delta}p)}(m-2,1)^2-4\eq0\ (\mo\ p).$$
By [S10a, Lemma 2.3], $v_{p-(\f{\Delta}p)}(m-2,1)\eq2\ (\mo\ p)$. So
$$v_{p-(\f{\Delta}p)}(m-2,1)\eq 2\ (\mo\ p^2).$$
By induction,
$(m-2)u_n(m-2,1)\pm v_n(m-2,1)=2u_{n\pm1}(m-2,1)$ for all $n\in\Z^+$. Therefore
$$\align 2u_p(m-2,1)=&(m-2)u_{p-(\f{\Delta}p)}(m-2,1)+\l(\f{\Delta}p\r)v_{p-(\f{\Delta}p)}(m-2,1)
\\\eq&(m-2)u_{p-(\f{\Delta}p)}(m-2,1)+2\l(\f{\Delta}p\r)\ (\mo\ p)
\endalign$$
and hence (2.6) follows.

The proof of Lemma 2.4 is now complete. \qed

Combining Lemmas 2.3 and 2.4 we get the following result.

\proclaim{Lemma 2.5} Let $p$ be an odd prime and let $m\in\Z$ with $\Delta=m(m-4)\not\eq0\ (\mo\ p)$. Then
$$\sum_{k=1}^{p-1}(-1)^k\f{\bi{2k}k}{k(m-4)^k}
\eq q_p(m-4)-\f m2\l(\f{\Delta}p\r)\f{u_{p-(\f{\Delta}p)}(m-2,1)}p\ (\mo\ p).\tag2.7$$
\endproclaim

\heading{3. Proofs of Theorems 1.1--1.2 and (1.3)}\endheading

\medskip
\noindent{\it Proof of Theorem 1.2}.
For $k=0,\ldots,p^a-1$, clearly
$$\align \binom{hp^a-1}k(-1)^k=&(-1)^k\prod_{0<j\ls k}\f{hp^a-j}j
=\prod_{0<j\ls k}\l(1-h\f{p^a}j\r)
\\\eq&1-h\sum_{0<j\ls k}\f{p^a}j=1-hp^aH_k\ (\mo\ p).
\endalign$$
Thus
$$\sum_{k=0}^{p^a-1}\bi{hp^a-1}k\f{\bi{2k} k}{(-m)^k}\eq\sum_{k=0}^{p^a-1}\f{\bi{2k} k}{m^k}
-hp^a\sum_{k=0}^{p^a-1}\f{H_k}{m^k}\bi{2k}k\ (\mo\ p).$$
If $p\nmid m-4$, then by applying Lemmas 2.1, 2.2 and 2.5 we get
$$\align&\sum_{k=0}^{p^a-1}\bi{hp^a-1}k\f{\bi{2k} k}{(-m)^k}
\\\eq&\sum_{k=0}^{p^a-1}\f{\bi{2k} k}{m^k}
+ph\l(\f{\Delta}{p^a}\r)\sum_{k=1}^{p-1}(-1)^k\f{\bi{2k}k}{k(m-4)^k}
\\\eq&\l(\f{\Delta}{p^a}\r)+\l(\f{\Delta}{p^{a-1}}\r)u_{p-(\f{\Delta}p)}(m-2,1)
\\&+ph\(\l(\f{\Delta}{p^a}\r)q_p(m-4)-\f m2\l(\f{\Delta}{p^{a-1}}\r)\f{u_{p-(\f{\Delta}p)}(m-2,1)}p\)
\ (\mo\ p)
\endalign$$
and hence (1.7) follows.
In the case $m\eq 4\ (\mo\ p)$, we have
$$p^a\sum_{k=1}^{p^a-1}\f{H_k}{m^k}\bi{2k}k\eq 2p\da_{a,1}\ (\mo\ p^2)$$
by Lemma 2.1, and
$$\align &u_{p-(\f{\Delta}p}(m-2,1)=u_p(m-2,1)
\\\eq& p\l(\f{m-2}2\r)^{p-1}+\da_{p,3}m\f{m-4}3
\eq p+\da_{p,3}(m-4)\ (\mo\ p^2).
\endalign$$
by [S10b, Lemma 2.2]. So (1.7) also holds when $p\mid m-4$.

Since $u_{p-(\f{\Delta}p)}(m-2,1)\eq0\ (\mo\ p)$ by [S10a, Lemma 2.3], (1.7) in the case $hm\eq2\ (\mo\ p)$ yields (1.8).

So far we have completed the proof of Theorem 1.2. \qed

\medskip
\noindent{\it Proof of Theorem 1.1}. Choose $m\in\Z$ such that $hm\eq2\ (\mo\ p^2)$.
Clearly $p\nmid m$. Note that
$$m-4\eq\f 2h-4=\f{2-4h}h\ (\mo\ p^2).$$
So we may get (1.1) by applying (1.8). This concludes the proof of Theorem 1.1. \qed

\medskip
\noindent{\it Proof of (1.3)}. For $k\in\N$ clearly the constant term of
$$(2-x-x^{-1})^k=\f{(-1)^k}{x^k}(x-1)^{2k}$$
is the central binomial coefficient $\bi{2k}k$.
Observe that
$$\sum_{k=0}^{n-1}\bi{n-1}k(-1)^km^{n-1-k}(2-x-x^{-1})^k=(m-2+x+x^{-1})^{n-1}.$$
Comparing the constant terms of both sides of the last equality we obtain
$$\align&\sum_{k=0}^{n-1}\bi{n-1}k\bi{2k}k(-1)^km^{n-1-k}
\\=&\sum_{k=0}^{\lfloor(n-1)/2\rfloor}\bi{n-1}{k,k,n-1-2k}(m-2)^{n-1-2k},
\endalign$$
which is equivalent to (1.3). We are done. \qed

\enddocument